\documentclass[a4paper,12pt]{amsart}
\usepackage[latin1]{inputenc}
\usepackage[T1]{fontenc}
\usepackage{amsfonts,amssymb}
\usepackage{pstricks}
\usepackage{amsthm}
\usepackage[all]{xy}
\usepackage[pdftex]{graphicx}
\usepackage{graphicx}   
\newtheorem{theo}{Theorem}[section]

\newtheorem{rem}[theo]{Remark}
\newtheorem{ex}{Example}

\title[Even-Odd partition identities of the Rogers-Ramanujan type]{Even-Odd partition identities of Rogers-Ramanujan type}

\author{Pooneh Afsharijoo}

\begin {document}
\maketitle
\begin{abstract} We prove a theorem which adds a new member to the \textit{Rogers-Ramanujan identities}. This new member counts partitions with different type of constraints on even and odd parts. Generalizing this theorem, we obtain two families of partition identities.
\end{abstract}
\section{Introduction}

A \textit{partition} $\Lambda$ of length $\ell(\Lambda)=m$ of a positive integer $n$ is a decreasing sequence of $m$ positive integers $(\lambda_1\geq \cdots \geq \lambda_m)$ whose sum is equal to $n.$ These positive integers are called the \textit{part} of the partition $\Lambda$. We denote by $p(n)$ the number of partitions of $n.$ By convention, zero has just one partition which is the empty partition.
%For example $4$ has the following five partitions (so $p(4)=5$):
%For example $5$ has the  following seven partitions (so $p(5)=7$):
%$$(4), \ (3,1),\  (2,2), \ (2,1,1), \ (1,1,1,1).$$

%$$(5), \ (4,1), \ (3,2), \ (3,1,1), \ (2,2,1), \ (2,1,1,1), \ (1,1,1,1,1).$$
   One important subject in partition theory is the study of the \textit{partition identities}. A partition identity is an equality between the number of partitions of an integer $n$ satisfying a property $P_1$ and the number of partitions of $n$ satisfying another property $P_2,$ which is true for every $n.$

 Two of the most famous partition identities are the \textit{Rogers-Ramanujan identities}:

\begin{theo}(Rogers-Ramanujan identities) 
For $i=1,2,$ let $B_{2,i}(n)$ denote the number of partitions of $n$ with no consecutive parts, neither equal parts and at most $i-1$ of parts are equal to $1$. Let $A_{2,i}(n)$ denote the number of partitions of $n$ into parts $\not\equiv 0,\pm i \ ( \text{mod}. 5)$. Then $A_{2,i}(n)=B_{2,i}(n)$ for all integers $n$.
\end{theo}

%\begin{ex}\label{ex}

%The set $\{(4,1), \ (1,1,1,1,1) \}$ is the set of all partitions of $5$ which are counted by $A_{2,2}(5)$ and the set $\{(5), \ (4,1)\}$ is the set of those partitions of $5$ which are counted by $B_{2,2}(5).$ So we have $A_{2,2}(5)=B_{2,2}(5)=2.$
%\end{ex}

In general, it is very difficult to \textit{guess and prove} partition identities. In \cite{AM}, we used the relation between generating series of partitions satisfying some conditions and the \textit{arc spaces} (this was established in \cite{BMS} and \cite{BMS1}; see also \cite{M} for relation with singularities), to guess and prove some new identities extending the Rogers-Ramanujan identities. This means that we have found some new type of partitions whose number is equal to $A_{2,i}(n)=B_{2,i}(n)$ in the Rogers-Ramanujan identities. We denote  the number of these partitions by $C_{2,i}(n).$
\\
\\
In \cite{Af}, we used a combinatorial method to prove this extension of the Rogers-Ramanujan identities. This last method was first given in \textit{ Andrews-Baxter system of recurrence relations}.  In \cite{A}, G. Andrews shows that the $B_{2,i}(n)$ are uniquely determined by some system of recurrence relations. In \cite{Af}, we prove that the $C_{2,i}(n)$ satisfy also this system in order to prove that $C_{2,i}(n)$ is equal to $B_{2,i}(n)$.
\\
\\
In this paper, using similar types of ideas we were able to create a new system of recurrence relations and find a new extension of the Rogers-Ramanujan identities in which even and odd parts play different roles.

To introduce this new member of the Rogers-Ramanujan identities we fix some notations. Let $n$ be an integer and $\Lambda$ be a partition of $n$. We regroup the even parts (respectively the odd parts) of $\Lambda$ together and we denote them by $\lambda_{i_1}\geq \lambda_{i_2} \geq \cdots \geq \lambda_{i_{r_1}}$ (respectively by $\lambda_{j_1}\geq \lambda_{j_2} \geq \cdots \geq \lambda_{j_{r_2}}$). 
We have the following  theorem (see Theorem \ref{k=1} below):

\begin{theo}
  Let  $P_{2,i}(n)$ denote the number of partitions of $n$ with at most $i-1$ parts equal to $1,$ whose smallest even part is greater than or equal to two times its length, and whose odd parts satisfy $\lambda_{j_\ell}-\lambda_{j_{\ell+2}}\geq 4.$ i.e.,
  $$P_{2,i}(n):=\{\Lambda:n \ | \ \lambda_{i_{r_1}}\geq 2\ell(\Lambda),\ \lambda_{j_\ell}-\lambda_{j_{\ell+2}}\geq 4, \text{at most $i-1$ parts equal to 1} \}.$$
  Then $P_{2,i}(n)=B_{2,i}(n)=A_{2,i}(n),$ where $A_{2,i}(n)$ and $B_{2,i}(n)$ are the same as in the Rogers-Ramanujan identities.
\end{theo}

\begin{ex} In this example we will show this theorem for the partitions of $6$ when $i=2.$ Note that $6$ has the following $11$ partitions: \begin{equation*}
\begin{aligned}
6 & =6 \\
      &=5+1\\
      &= 4+2\\
      &=4+1+1\\
      &=3+3\\
      &=3+2+1\\
      &=3+1+1+1\\
      &=2+2+2\\
      &=2+2+1+1\\
      &=2+1+1+1+1\\
      &=1+1+1+1+1+1.
\end{aligned}
\end{equation*}
\begin{itemize}
\item The partitions of $6$ which are counted by $A_{2,2}(6)$ are those with parts $\equiv 1,4  \ ( \text{mod}. 5).$ These are the following three partitions:
$$6, \ \ 4+1+1, \ \ 1+1+1+1+1+1.$$ 
\item The partitions of $6$ without equal or consecutive parts are the following three partitions:
$$6, \ \ 5+1,\ \ 4+2.$$
\item The partition of $6$ which are counted by $P_{2,2}(6)$ are the following three partitions:

$$6, \ \ 5+1, \ \ 3+3.$$
In the all other partitions of $6$ either the smallest even part is strictly less than 2 times the length, or there is at least a block of three odd parts with the difference between the first and the third $<4.$

\end{itemize}
Therefore we have $P_{2,2}(6)=A_{2,2}(6)=B_{2,2}(6)=3.$

%In Example (\ref{ex}) we saw that $A_{2,2}(5)=B_{2,2}(5)=2.$
 
%5 Note that $(5)$ and $(4,1)$ are the only partitions of $5$ which are counted by $P_{2,2}(5).$ So we have $P_{2,2}(5)=A_{2,2}(n)=B_{2,2}(5)=2.$
\end{ex}

 In the second section, we prove this theorem. To do so, we denote by $p_{2,i}(m,n)$ (respectively by $b_{2,i}(m,n)$) the number of partitions of $n$  which are counted by $P_{2,i}(n)$ (respectively by $B_{2,i}(n)$) with exactly $m$ parts. Then we construct a new system of  recurrence relations  between $p_{2,i}(m,n)$ and we prove that $b_{2,i}(m,n)$ satisfy the same system of equations (see Theorem \ref{k=1}). 
 
 In the last section, we generalize Theorem \ref{k=1} and first we obtain the following family of partition identities (see Theorem \ref{2k+1} below): 
\begin{theo} For all integers $k,n\geq 0$ and $i=1$ or $2$ let $P_{2,i}^{2k+1}(n)$ denote the number of partitions of $n$ whose parts are greater than or equal to $2k+1,$ with at most $i-1$ parts equal to $2k+1,$ whose smallest even part is greater than or equal to two times its length plus $k$, and whose odd parts satisfy $\lambda_{j_\ell}-\lambda_{j_{\ell+2}}\geq 4.$ i.e.,
 
 \begin{equation*}
 \begin{split}
	P_{2,i}^{2k+1}(n):=\{&\Lambda:n \ | \ \lambda_{i_{r_1}},\lambda_{j_{r_2}} \geq 2k+1,  \ \lambda_{i_{r_1}}\geq 2(\ell(\Lambda)+k),\\
	&\lambda_{j_\ell}-\lambda_{j_{\ell+2}}\geq 4, \text{at most $i-1$ parts equal to $2k+1$} \}.
\end{split}
\end{equation*}

  Let $B_{2,i}^{2k+1}(n)$ denote the number of partitions of $n$ whose parts are greater than or equal to $2k+1,$ with at most $i-1$ parts equal to $2k+1$ and without equal or consecutive parts. Then $B_{2,i}^{2k+1}(n)=P_{2,i}^{2k+1}(n).$ 
 \end{theo}

 %Note that if we take $k=0$ we obtain a new member of Rogers-Ramanujan identities.
In order to prove this theorem we define a simple bijective transformation between partitions to send each partition which is counted by $P_{2,i}^{2k+1}(n)$ (respectively by $B_{2,i}^{2k+1}(n)$) to a partition which is counted by $P_{2,i}(n)$ (respectively by $B_{2,i}(n)$). Then we apply Theorem \ref{k=1}. 
\\
\\
Finely, using Theorem \ref{2k+1} we give the following family of partitions identities (see Theorem \ref{2k} below):

 \begin{theo} For all integers $n\geq 0, k\geq 1$ and $i=1$ or $2$ let $P_{2,i}^{2k}(n)$ denote the number of partitions of $n$ whose parts are greater than or equal to $2k,$ with at most $i-1$ parts equal to $2k,$ whose smallest odd part plus $1$ is greater than or equal to two times its length plus $k$, and whose even parts satisfy $\lambda_{i_\ell}-\lambda_{i_{\ell+2}}\geq 4.$ i.e.,

 \begin{equation*}
 \begin{split}
	P_{2,i}^{2k}(n):=\{&\Lambda:n \ | \ \lambda_{i_{r_1}},\lambda_{j_{r_2}} \geq 2k,  \ \lambda_{j_{r_2}}+1\geq 2(\ell(\Lambda)+k),\\
	&\lambda_{i_\ell}-\lambda_{i_{\ell+2}}\geq 4, \text{at most $i-1$ parts equal to $2k$} \}.
\end{split}
\end{equation*} 
 
  Let $B_{2,i}^{2k}(n)$ denote the number of partitions of $n$ whose parts are greater than or equal to $2k,$ with at most $i-1$ parts equal to $2k$ and without equal or consecutive parts. Then $B_{2,i}^{2k}(n)=P_{2,i}^{2k}(n).$ 
 \end{theo}

\section*{ACKNOWLEDGMENT} I would like to express my special thanks to my Ph.D adviser, Hussein Mourtada, for suggesting me this project, his motivation and supports. His continues guidance helped me to write this paper. I also would like to thank Jehanne Dousse, Frederic Jouhet and Bernard Teissier with whom I had the chance to discuss about partition identities. 
\section{Even-odd new member of the Rogers-Ramanujan identities}
In this section we give a new member of the Rogers-Ramanujan identities whose behavior on even and odd parts of a partition is different. We use the notations used in the introduction: 
\begin{theo}\label{k=1}
  Let  $P_{2,i}(n)$ denote the number of partitions of $n$ with at most $i-1$ parts equal to $1,$ whose smallest even part is greater than or equal to two times its length, and whose odd parts satisfy $\lambda_{j_\ell}-\lambda_{j_{\ell+2}}\geq 4.$ i.e.,
  $$P_{2,i}(n):=\{\Lambda:n \ | \ \lambda_{i_{r_1}}\geq 2\ell(\Lambda),\ \lambda_{j_\ell}-\lambda_{j_{\ell+2}}\geq 4, \text{at most $i-1$ parts equal to 1} \}.$$
  Then $P_{2,i}(n)=B_{2,i}(n),$ where $B_{2,i}(n)$ is the same as in the Rogers-Ramanujan identities.
\end{theo}
 
 \begin{proof}
 Let $p_{2,i}(m,n)$ (respectively $b_{2,i}(m,n)$) denote the number of partitions of $n$ with exactly $m$ parts and which are counted by $P_{2,i}(n)$ (respectively by $B_{2,i}(n)$). We prove that $p_{2,i}(m,n)$ and $b_{2,i}(m,n)$ both satisfy the following system:
 
 \begin{equation}\label{1}
 \begin{split}
 &p_{2,i}(m,n)=\begin{cases}
1 &\text{ if } m=n=0  \\
0 &\text{ if } m\leq 0 \text{ or } n\leq 0 \text{ but } (m,n) \neq (0,0);
\end{cases}\\
&p_{2,2}(m,n)-p_{2,1}(m,n)=p_{2,2}(m-1,n-2m+1); \\
&p_{2,1}(m,n) =p_{2,1}(m-1,n-2m)+p_{2,2}(m,n-2m).
\end{split}
\end{equation} 

Note that $0$ has only one partition whose length is zero (the empty partition). A negative number has no partition, and a positive number has no partition of non positive length. So the first equation is true for $b_{2,i}(m,n)$ and $p_{2,i}(m,n).$
\\
\\
In order to prove the second equation of System (\ref{1}) for  $p_{2,i}(m,n)$  we define a bijection between the partitions counted by each side of this equation. Note that the left hand side of this equation counts the number of partitions $\Lambda$ of $n$ with exactly $m$ parts and exactly one part equal to $1,$ whose smallest even part $\geq 2m$, and whose odd parts satisfy $\lambda_{j_\ell}-\lambda_{j_{\ell+2}}\geq 4.$ i.e.,   
 $$\Lambda:(\underbrace{\lambda_{i_1}\geq \lambda_{i_2} \geq \cdots \geq \lambda_{i_{r_1}}}_{\text {The even parts of $\Lambda$}},\underbrace{\lambda_{j_1}\geq \lambda_{j_2} \geq \cdots \geq \lambda_{j_{r_2-1}}\geq 1}_{\text {The odd parts of $\Lambda$}}),$$
 where $r_1+r_2=m, \ \lambda_{i_{r_1}}\geq 2m, \ \lambda_{j_{r_{2}-1}}\geq3$ and $\lambda_{j_\ell}-\lambda_{j_{\ell+2}}\geq 4.$
We transform $\Lambda$ to a partition $\mu$ by deleting its smallest part (which is equal to one) and removing $2$ from all other parts. i.e.,

$$\mu:(\underbrace{\lambda_{i_1}-2\geq \lambda_{i_2}-2 \geq \cdots \geq \lambda_{i_{r_1}}-2}_{\text {The even parts of $\mu$}},\underbrace{\lambda_{j_1}-2\geq \lambda_{j_2}-2 \geq \cdots \geq \lambda_{j_{r_2-1}}-2}_{\text {The odd parts of $\mu$}}).$$
 Thus, if $m=1$, then $\Lambda:(1)$ and so $\mu$ is the empty partition; If $m\geq2 $ we obtain a partition of $n-2m+1$ with exactly $m-1$ parts, whose smallest even part is equal to $\lambda_{i_{r_1}}-2$ which is greater than or equal to $2m-2=2(m-1);$ its odd parts are greater than or equal to $1$ and satisfy:
$$(\lambda_{j_\ell}-2)-(\lambda_{j_{\ell+2}}-2)=\lambda_{j_\ell}-\lambda_{j_{\ell+2}}\geq 4.$$
Note that if $\lambda_{j_{r_2-1}}=3,$ since $\lambda_{j_{r_2}}=1$ so $\lambda_{j_{r_2-2}} \neq 3.$ This means that $\mu$ has at most one part equal to $1.$ So $\mu$ is a partition which is counted by $p_{2,2}(m-1,n-2m+1).$ Obviously this transformation is a bijection which proves the second equation of System (\ref{1}) for $p_{2,i}(m,n).$
\\
\\
In order to prove the last equation we take a partition $$\Lambda:(\underbrace{\lambda_{i_1}\geq \lambda_{i_2} \geq \cdots \geq \lambda_{i_{r_1}}}_{\text {The even parts of $\Lambda$}},\underbrace{\lambda_{j_1}\geq \lambda_{j_2} \geq \cdots \geq \lambda_{j_{r_2}}}_{\text {The odd parts of $\Lambda$}}),$$ which is counted by $p_{2,1}(m,n).$

\begin{itemize}
\item If $\lambda_{i_{r_1}}=2m:$ then we transform $\Lambda$ to a partition $\mu_1$ by deleting $\lambda_{i_{r_1}}$ from $\Lambda.$ We obtain a partition of $n-2m$ with exactly $m-1$ parts whose smallest even part (respectively odd part) is equal to $\lambda_{i_{r_1-1}}$ (respectively $\lambda_{j_{r_2}}$) and we have:
$$\lambda_{i_{r_1-1}}\geq \lambda_{i_{r_1}}=2m>2(m-1).$$
So $\mu_1$ is a partition which is counted by $p_{2,1}(m-1,n-2m)$ whose smallest even part is strictly greater than its length.
\\
\item If $\lambda_{i_{r_1}} \geq 2m+2$ and $\lambda_{j_{r_2}}=\lambda_{j_{r_2-1}}=3:$ then we transform $\Lambda$ to a partition $\mu_2$ by deleting its last two odd parts, adding a smallest even part equal to $2m-2$ and removing $4$ from all other parts. i.e.,

% \begin{equation*}
% \begin{split}
	%\mu_2:&(\underbrace{\lambda_{i_1}-4\geq \lambda_{i_2}-4 \geq \cdots \geq \lambda_{i_{r_1}}-4 \geq 2m-2}_{\text {The even parts of $\mu$}}\\
	%&\underbrace{\lambda_{j_1}-4\geq \lambda_{j_2}-4 \geq \cdots \geq \lambda_{j_{r_2-2}}-4}_{\text {The odd parts of $\mu$}}).
%\end{split}
%\end{equation*} 

$$\mu_2:(\underbrace{\lambda_{i_1}-4\geq \lambda_{i_2}-4 \geq \cdots \geq \lambda_{i_{r_1}}-4 \geq 2m-2}_{\text {The even parts of $\mu_2$}},\underbrace{\lambda_{j_1}-4\geq \lambda_{j_2}-4 \geq \cdots \geq \lambda_{j_{r_2-2}}-4}_{\text {The odd parts of $\mu_2$}}).$$ 

We obtain a partition of $n-2m$ with exactly $m$ parts whose smallest even part (respectively odd part) is equal to $2m-2$ which is equal to two times the length of $\mu_2$ (respectively is equal to $\lambda_{j_{r_2-2}}-4$ which is greater than or equal to $3$). We have also: 
$$(\lambda_{j_\ell}-4)-(\lambda_{j_{\ell+2}}-4)=\lambda_{j_\ell}-\lambda_{j_{\ell+2}}\geq 4.$$
So in this case $\mu_2$ is a partition which is counted by $p_{2,1}(m-1,n-2m)$ whose smallest even part is equal to its length.
\\
\item If $\lambda_{i_{r_1}} \geq 2m+2$ and $\Lambda$ has at most one part equal to $3$: then we transform $\Lambda$ to a partition $\mu_3$ by removing $2$ from each part. We obtain a partition of $n-2m$ with exactly $m$ parts whose smallest even part $\lambda_{i_{r_1}}-2$ is greater than or equal to $2m$, whose smallest odd part $\lambda_{j_{r_2}}-2$ is greater than or equal to $1,$ with at most one part equal to $1$ and whose odd parts satisfy the following inequality:
$$(\lambda_{j_\ell}-2)-(\lambda_{j_{\ell+2}}-2)=\lambda_{j_\ell}-\lambda_{j_{\ell+2}}\geq 4.$$
So $\mu_3$ is a partition which is counted by $p_{2,2}(m,n-2m).$
\end{itemize}
Obviously the last three transformations defined above are bijective and they prove the last equation of System (\ref{1}) for $p_{2,i}(m,n).$ 
\\
\\
So far we proved that the $p_{2,i}(m,n)$ satisfy System (\ref{1}). We prove it now for the $b_{2,i}(m,n).$ Let $\Lambda$ be a partition which is counted by $b_{2,2}(m,n)-b_{2,1}(m,n).$ So it is a partition of $n$ with exactly $m$ parts and one part equal to $1,$ without equal or consecutive parts. We send $\Lambda$ to a partition $\mu$ by removing $2$ from each part. We obtain a partition which is counted by $b_{2,2}(m-1,n-m).$ This transformation defines a bijection between the partition counted by each side of the second equation of System (\ref{1}) for $b_{2,i}(m,n).$ 
\\
\\
In order to prove the last equation of this system for $b_{2,i}(m,n)$, we take a partition $\Lambda:(\lambda_1\geq \cdots \geq \lambda_m)$ which is counted by $b_{2,1}(m,n).$ So it is a partition of $n$  without consecutive or equal parts and whose parts are greater than or equal to $2.$ We transform it to a partition $\mu$ by removing $2$ from each part.
\begin{itemize}
\item If $\lambda_m=2:$ then we obtain a partition of $n-2m$ with exactly $m-1$ parts $\geq 2,$ without consecutive or equal parts. So $\mu$ is counted by $b_{2,1}(m-1,n-2m).$ 
\\
\item If $\lambda_m\geq 3:$ then we obtain a partition of $n-2m$ with exactly $m$ parts $\geq 1,$ without equal or consecutive parts. So in this case $\mu$ is counted by $b_{2,2}(m,n-2m).$
\end{itemize}
Note that the last two deformations defined above are bijective and they prove that $b_{2,i}(m,n)$ satisfy the last equation of System (\ref{1}).
\\
So far we proved that $b_{2,i}(m,n)$ and $p_{2,i}(m,n)$ both satisfy System (\ref{1}). By double induction on $n,m$ one can show that the $p_{2,i}(m,n)$ are \textit{uniquely} determined by System (\ref{1}). Therefore, $p_{2,i}(m,n)=b_{2,i}(m,n)$ for all integers $m,n$ and $i=1$ or $2.$ So we have:

$$P_{2,i}(n)=\sum_{m\geq 0} p_{2,i}(m,n)=\sum_{m\geq0} b_{2,i}(m,n)=B_{2,i}(n).$$

 \end{proof}
 
 \begin{rem} Note that by proving $p_{2,i}(m,n)=b_{2,i}(m,n)$ for all integers $m,n$ and $i=1$ or $2,$ actually we proved that even if we fix the length of the partitions of $n$, the equality between $P_{2,i}(n)$ and $B_{2,i}(n)$ holds. This is not true in general for $A_{2,i}(n)$ and $B_{2,i}(n).$ 
 \end{rem}

 %By generalizing System (\ref{1}) in the proof of Theorem (\ref{k=1}) we obtain a family of partition identities of Rogers-Ramanujan type:
 \section{Generalization of Theorem \ref{k=1}}
 In This section we give two families of partition identities by generalizing Theorem \ref{k=1}. The first one is as follows:
 \begin{theo}\label{2k+1} For all integers $n\geq 0, \ k\geq 1$ and $i=1$ or $2$ let $P_{2,i}^{2k+1}(n)$ denote the number of partitions of $n$ whose parts are greater than or equal to $2k+1,$ with at most $i-1$ parts equal to $2k+1,$ whose smallest even part is greater than or equal to two times its length plus $k$, and whose odd parts satisfy $\lambda_{j_\ell}-\lambda_{j_{\ell+2}}\geq 4.$ i.e.,
\begin{equation*}
\begin{split} 
	P_{2,i}^{2k+1}(n):=\{&\Lambda:n \ | \ \lambda_{i_{r_1}},\lambda_{j_{r_2}} \geq 2k+1,  \ \lambda_{i_{r_1}}\geq 2(\ell(\Lambda)+k),\\ 
	&\lambda_{j_\ell}-\lambda_{j_{\ell+2}}\geq 4, \text{at most $i-1$ parts equal to $2k+1$} \}.
\end{split}
\end{equation*}

  Let $B_{2,i}^{2k+1}(n)$ denote the number of partitions of $n$ whose parts are greater than or equal to $2k+1,$ with at most $i-1$ parts equal to $2k+1$ and without equal or consecutive parts. Then $B_{2,i}^{2k+1}(n)=P_{2,i}^{2k+1}(n).$ 
 \end{theo}

 \begin{proof} Let denote by $p_{2,i}^{2k+1}(m,n)$ (respectively by $b_{2,i}^{2k+1}(m,n)$) the number of partitions of $n$ which are counted by $P_{2,i}^{2k+1}(n)$ (respectively by $B_{2,i}^{2k+1}(n)$) with exactly $m$ parts. We prove that $p_{2,i}^{2k+1}(m,n)$ and $b_{2,i}^{2k+1}(m,n)$ satisfy the following equations for all $k\geq 1$:
 
 \begin{equation}\label{p1}
 p_{2,i}^{2k+1}(m,n)=p_{2,i}(m,n-2mk),
\end{equation} 
 and 
 
 \begin{equation}\label{b1}
  b_{2,i}^{2k+1}(m,n)=b_{2,i}(m,n-2mk).
\end{equation}  

In order to prove Equation (\ref{p1}) we define a transformation $T$ from the set of all partitions which are counted by $ p_{2,i}^{2k+1}(m,n)$ to the set of all partitions which are counted by $p_{2,i}(m,n-2mk)$ as follows:  
\\
Let $$\Lambda:(\underbrace{\lambda_{i_1}\geq \lambda_{i_2} \geq \cdots \geq \lambda_{i_{r_1}}}_{\text {The even parts of $\Lambda$}},\underbrace{\lambda_{j_1}\geq \lambda_{j_2} \geq \cdots \geq \lambda_{j_{r_2}}}_{\text {The odd parts of $\Lambda$}}),$$ be a partition which is counted by $ p_{2,i}^{2k+1}(m,n).$ So it has exactly $m$ parts, each greater than or equal to $2k+1,$ with at most $i-1$ part equal to $2k+1$, whose smallest even part $\lambda_{i_{r_1}}\geq 2(m+k)$ and whose odd parts satisfy the following inequality:
$$\lambda_{j_\ell}-\lambda_{j_{\ell+2}}\geq 4.$$
We remove $2k$ from each part of $\Lambda$ and we obtain a partition $\mu$ as follows:

$$T(\Lambda):=\mu:(\underbrace{\lambda_{i_1}-2k\geq \lambda_{i_2}-2k \geq \cdots \geq \lambda_{i_{r_1}}-2k}_{\text {The even parts of $\mu$}},\underbrace{\lambda_{j_1}-2k\geq \lambda_{j_2}-2k \geq \cdots \geq \lambda_{j_{r_2}}-2k}_{\text {The odd parts of $\mu$}}),$$

with exactly $m$ parts (each $\geq 1$), whose smallest even part $\mu_{i_{r_1}}=\lambda_{i_{r_1}}-2k$ is greater than or equal to $2m.$  Moreover, its odd parts satisfy the following inequality:

$$\mu_{j_\ell}-\mu_{j_{\ell+2}} =(\lambda_{j_\ell}-2k)-(\lambda_{j_{\ell+2}}-2k) = \lambda_{j_\ell}-\lambda_{j_{\ell+2}}\geq 4.$$

So $\mu$ is a partition which is counted by $p_{2,i}(m,n-2mk)$. The transformation $T$ is obviously injective. In order to show that it is also surjective we take a partition

% Note that if we add $2k$ to each part of $\mu,$ we obtain the partition $\Lambda.$ so this transformation is a bijection. This proves Equation (\ref{p1}).
%We define now the inverse of this transformation. Let  
$$\mu:(\underbrace{\mu_{i_1}\geq \mu_{i_2} \geq \cdots \geq \mu_{i_{r_1}}}_{\text {The even parts of $\mu$}},\underbrace{\mu_{j_1}\geq \mu_{j_2} \geq \cdots \geq \mu_{j_{r_2}}}_{\text {The odd parts of $\mu$}}),$$ which is counted by $ p_{2,i}(m,n-2mk).$ So it has exactly $m$ parts, (each greater than or equal to $1$), with at most $i-1$ part equal to $1$, whose smallest even part $\mu_{i_{r_1}}\geq 2m$ and whose odd parts satisfy the following inequality:
$$\mu_{j_\ell}-\mu_{j_{\ell+2}}\geq 4.$$
We add $2k$ to each part of $\mu$ and we obtain a partition $\Lambda$ as follows:
$$\Lambda:(\underbrace{\mu_{i_1}+2k\geq \mu_{i_2}+2k \geq \cdots \geq \mu_{i_{r_1}}+2k}_{\text {The even parts of $\Lambda$}},\underbrace{\mu_{j_1}+2k\geq \mu_{j_2}+2k \geq \cdots \geq \mu_{j_{r_2}}+2k}_{\text {The odd parts of $\Lambda$}}),$$
with exactly $m$ parts, each $\geq 2k+1$, whose smallest even part $\lambda_{i_{r_1}}=\mu_{i_{r_1}}+2k$ is greater than or equal to $2(m+k).$  Moreover, its odd parts satisfy the following inequality:
$$\lambda_{j_\ell}-\lambda_{j_{\ell+2}} =(\mu_{j_\ell}+2k)-(\mu_{j_{\ell+2}}+2k) = \mu_{j_\ell}-\mu_{j_{\ell+2}}\geq 4.$$
Thus, $\Lambda$ is a partition which is counted by $ p_{2,i}^{2k+1}(m,n)$ and we have $T(\Lambda)=\mu.$ So the transformation $T$ is a bijection. This proves Equation (\ref{p1}).
%and the bijectivity of this transformation from $\Lambda$ to $\mu$ proves Equation (\ref{p1}).
\\
\\
Now let $\Lambda:(\lambda_1 \geq \cdots \geq \lambda_m)$ be a partition which is counted by the left hand side of Equation (\ref{b1}). So it has not equal or consecutive parts, with  $\lambda_m\geq 2k+1$ and at most $i-1$ parts equal to $2k+1.$ We send $\Lambda$ to a partition $\mu$ using again the transformation $T.$ We obtain a partition of $n-2mk$ with exactly $m$ parts, at most $(i-1)$ part equal to $1$ and without equal or consecutive parts. So $\mu$ is a partition which is counted by $b_{2,i}(m,n-2mk).$ Obviously the transformation $T$ defines a one to one correspondence between $\Lambda$ and $\mu$ and so proves Equation (\ref{b1}). Now we have:

\begin{align*}
 p_{2,i}^{2k+1}(m,n) \underset{\text{By Equation (\ref{p1})} }{=}& p_{2,i}(m,n-2mk)\\
 \\
  \underset{\text{By Theorem \ref{k=1}} }{=}& b_{2,i}(m,n-2mk)\\
  \\
   \underset{\text{By Equation  (\ref{b1})} }{=}& b_{2,i}^{2k+1}(m,n).\\
\end{align*}

This last equation gives us:

$$P_{2,i}^{2k+1}(n)=\sum_{m\geq 0} p_{2,i}^{2k+1}(m,n)=\sum_{m\geq0} b_{2,i}^{2k+1}(m,n)=B_{2,i}^{2k+1}(n).$$

 \end{proof}
 
 \begin{rem} Note that we can obtain the following system of equations between $p_{2,i}^{2k+1}(m,n)$ by generalizing System (\ref{1}):
 
  \begin{equation}\label{2}
 \begin{split}
 &p_{2,i}^{2k+1}(m,n)=\begin{cases}
1 &\text{ if } m=n=0  \\
0 &\text{ if } m\leq 0 \text{ or } n\leq 0 \text{ but } (m,n) \neq (0,0);
\end{cases}\\
\\
&p_{2,2}^{2k+1}(m,n)-p_{2,1}^{2k+1}(m,n)=p_{2,2}^{2k+1}(m-1,n-2m-2k+1); \\
\\
&p_{2,1}^{2k+1}(m,n) =p_{2,1}^{2k+1}(m-1,n-2m-2k)+p_{2,2}^{2k+1}(m,n-2m).
\end{split}
\end{equation} 
 Therefore, another proof of Theorem \ref{2k+1} is to show that $p_{2,i}^{2k+1}(m,n)$  and $b_{2,i}^{2k+1}(m,n),$ both satisfy the system above  by defining the similar transformations between partitions as in the proof of Theorem \ref{k=1}.
 \end{rem}
 Using this result, we prove another family of partition identities as follows:

 \begin{theo}\label{2k} For all integers $n\geq 0, k\geq 1$ and $i=1$ or $2$ let $P_{2,i}^{2k}(n)$ denote the number of partitions of $n$ whose parts are greater than or equal to $2k,$ with at most $i-1$ parts equal to $2k,$ whose smallest odd part plus $1$ is greater than or equal to two times its length plus $k$, and whose even parts satisfy $\lambda_{i_\ell}-\lambda_{i_{\ell+2}}\geq 4.$ i.e.,
\begin{equation*}
\begin{split}
P_{2,i}^{2k}(n):=\{&\Lambda:n \ | \ \lambda_{i_{r_1}},\lambda_{j_{r_2}} \geq 2k,  \ \lambda_{j_{r_2}}+1\geq 2(\ell(\Lambda)+k),\\ 
	&\lambda_{i_\ell}-\lambda_{i_{\ell+2}}\geq 4, \text{at most $i-1$ parts equal to $2k$} \}.
\end{split}
\end{equation*}
  Let $B_{2,i}^{2k}(n)$ denote the number of partitions of $n$ whose parts are greater than or equal to $2k,$ with at most $i-1$ parts equal to $2k$ and without equal or consecutive parts. Then $B_{2,i}^{2k}(n)=P_{2,i}^{2k}(n).$

 \end{theo}
 \begin{proof}
Let denote by $p_{2,i}^{2k}(m,n)$ (respectively by $b_{2,i}^{2k}(m,n)$) the number of partitions of $n$ which are counted by $P_{2,i}^{2k}(m,n)$ (respectively by $B_{2,i}^{2k}(m,n)$) and with exactly $m$ parts. We show first that for all $k\geq1$ we have: 
\begin{equation}\label{p}
 p_{2,i}^{2k}(m,n)=p_{2,i}^{2k+1}(m,n+m),
\end{equation} 
 and 
 
 \begin{equation}\label{b}
 b_{2,i}^{2k}(m,n)=b_{2,i}^{2k+1}(m,n+m).
\end{equation}  
In order to show Equation (\ref{p}) let 
 $$\Lambda:(\underbrace{\lambda_{i_1}\geq \lambda_{i_2} \geq \cdots \geq \lambda_{i_{r_1}}}_{\text {The even parts of $\Lambda$}},\underbrace{\lambda_{j_1}\geq \lambda_{j_2} \geq \cdots \geq \lambda_{j_{r_2}}}_{\text {The odd parts of $\Lambda$}}),$$ be a partition which is counted by  $p_{2,i}^{2k}(m,n)$. So its parts $\geq 2k,$ at most $i-1$ of them is equal to $2k$ and we have:
 
 $$r_1+r_2=m, \  \lambda_{i_\ell}-\lambda_{i_{\ell+2}}\geq 2, \ \text{and} \ \lambda_{j_{r_2}}+1\geq 2(m+k). $$
 
 We add $1$ to each part of $\Lambda$ and we obtain a partition of $n+m$ with exactly $m$ parts as follows:
 
  $$\mu:(\underbrace{\lambda_{i_1}+1\geq \lambda_{i_2}+1 \geq \cdots \geq \lambda_{i_{r_1}}+1}_{\text {The odd parts of $\mu$}},\underbrace{\lambda_{j_1}+1\geq \lambda_{j_2}+1 \geq \cdots \geq \lambda_{j_{r_2}}+1}_{\text {The even parts of $\mu$}}).$$ 
  
Note that all parts of $\mu$ are greater than or equal to $2k+1,$ its smallest even part $\geq 2(m+k)$ and the difference between the first and the third part of each block of three odd parts is greater than or equal to $4.$ So $\mu$ is a partition which is counted by $p_{2,i}^{2k+1}(m,n+m).$ Since this transformation from $\lambda$ to $\mu$ is obviously a bijection, we have Equation (\ref{p}).
\\
With the same transformation as below, we can send each partition $\Lambda$ which is counted by $b_{2,i}^{2k}(m,n)$ to a partition $\mu$ which is counted by $b_{2,i}^{2k+1}(m,n+m).$ Once again this  transformation defines a bijection and proves Equation (\ref{b}). So for all $k\geq1$ we have:
\\

\begin{align*}
 p_{2,i}^{2k}(m,n) \underset{\text{By Equation (\ref{p})} }{=}& p_{2,i}^{2k+1}(m,n+m)\\
 \\
  \underset{\text{By Theorem  \ref{2k+1}} }{=}& b_{2,i}^{2k+1}(m,n+m)\\
  \\
   \underset{\text{By Equation  (\ref{b})} }{=}& b_{2,i}^{2k}(m,n).\\
\end{align*}

Which gives us the following equations:

$$P_{2,i}^{2k}(n)=\sum_{m\geq 0} p_{2,i}^{2k}(m,n)=\sum_{m\geq0} b_{2,i}^{2k}(m,n)=B_{2,i}^{2k}(n).$$

 \end{proof}
 
  \begin{rem} Note that we can obtain the following system of equations between $p_{2,i}^{2k}(m,n)$ by generalizing System (\ref{1}):
 
  \begin{equation}\label{2}
 \begin{split}
 &p_{2,i}^{2k}(m,n)=\begin{cases}
1 &\text{ if } m=n=0  \\
0 &\text{ if } m\leq 0 \text{ or } n\leq 0 \text{ but } (m,n) \neq (0,0);
\end{cases}\\
\\
&p_{2,2}^{2k}(m,n)-p_{2,1}^{2k}(m,n)=p_{2,2}^{2k}(m-1,n-2m-2k+2); \\
\\
&p_{2,1}^{2k}(m,n) =p_{2,1}^{2k}(m-1,n-2m-2k+1)+p_{2,2}^{2k}(m,n-2m).
\end{split}
\end{equation} 
 Therefore, another proof of Theorem \ref{2k} is to show that $p_{2,i}^{2k}(m,n)$  and $b_{2,i}^{2k}(m,n),$ both satisfy the system above  by defining the similar transformations between partitions as in the proof of Theorem \ref{k=1}.
 \end{rem}

\newpage
\bibliographystyle{acm}
\nocite{AB, ADJM, G, GP, LZ}
\bibliography{article}

\end{document}